\documentclass[11pt]{amsart}

\usepackage{color}

\begin{document}

\numberwithin{equation}{section}

\def\R{{\mathbb R}}
\def\K{{\mathcal K}}
\def\L{{\mathcal L}}
\def\O{{\Omega}}
\def\I{{\mathcal I}}
\def\T{{\mathcal T}}
\def\vep{{\varepsilon}}
\def\p{{\partial}}
\def\a{{\alpha}}
\def\b{{\beta}}
\def\e{{\eta}}
\def\g{{\gamma}}
\def\s{{\sigma}}
\def\t{{\tau}}
\def\E{{\mathbb E}}
\def\N{{\mathbb N}}
\def\F{{\mathcal F}}
\def\P{{\mathbb P}}
\def\EE{{\mathcal E}}
\def\o{{\omega}}
\def\l{{\ell}}
\def\vp{{\varphi}}
\def\un{{\mathrm{1~\hspace{-1.4ex}l}}}
\def\to{{\rightarrow}}
\def\th{\theta}
\def\hs{H^s}
\def\Sd{S_{\Delta t}}
\def\Td{T_{\Delta t}}
\def\D{\Delta t}
\def\ds{\displaystyle}
\def\ud{u^{\D}}
\def\udR{u^{\D,R}}
\def\be{\begin{equation}}
\def\ee{\end{equation}}
\def\ba{\begin{array}}
\def\ea{\end{array}}
\def\l{\left}
\def\r{\right}

\newcommand{\Frac}{\displaystyle \frac}
\newcommand{\Sum}{\displaystyle \sum}
\newcommand{\Int}{\displaystyle \int}
\newcommand{\Sup}{\displaystyle \sup}

\newtheorem{Theorem}{Theorem}[section]
\newtheorem{Definition}[Theorem]{Definition}
\newtheorem{Proposition}[Theorem]{Proposition}
\newtheorem{Lemma}[Theorem]{Lemma}
\newtheorem{Corollary}[Theorem]{Corollary}
\newtheorem{Remark}[Theorem]{Remark}
\newtheorem{Example}[Theorem]{Example}
\newtheorem{Hypothesis}[Theorem]{Hypothesis}

\title[White noise dispersion for critical NLS]
{1D quintic nonlinear Schrodinger equation with white noise dispersion }

\author[A. Debussche]
{Arnaud Debussche } 
\address{IRMAR et ENS de Cachan, Antenne de Bretagne,
Campus de Ker Lann, Av. R. Schuman,
35170 BRUZ, FRANCE}
\email{arnaud.debussche@bretagne.ens-cachan.fr}
\author[Y. Tsutsumi]
{Yoshio Tsutsumi}
\address{Department of Mathematics, Kyoto University, Kyoto 606-8502, JAPAN}
\email{tsutsumi@math.kyoto-u.ac.jp}



\begin{abstract} In this article, we improve the Strichartz estimates obtained in \cite{dBD10} for 
the Schr\"odinger equation with white noise dispersion 
in one dimension. This allows us to prove global well posedness 
when a quintic critical nonlinearity is added to the equation. We finally show that the white noise
dispersion is the limit of smooth random dispersion. 
\end{abstract}

\maketitle


\section{Introduction}

The nonlinear Schr\"odinger equation with power nonlinearity is a common model in optics. It describes the propagation of waves in a nonlinear dispersive medium. It has been 
widely studied (see for instance \cite{cazenave}, \cite{sulem-sulem}). In the case
of a focusing nonlinearity, it has the form 
$$\left\{
\begin{array}{l}
\ds i \frac{du}{dt} +\Delta u  + |u|^{2\s}u =0,\; x\in\R,\; t>0,\\
\\
u(0)=u_0, \; x\in\R^n.
\end{array}
\right.
$$
It is well known that for subcritical nonlinearity, {\it i.e.} $\s<2/n$, this equation is globally well
posed in $L^2(\R^n)$ and in $H^1(\R^n)$ (\cite{Gi-Ve}, \cite{kato}, \cite{tsutsumi}). Moreover,
solitary waves are stable. 

For critical, $\s=2/n$, or supercritical, $\s>2/n$, nonlinearity, the equation is locally well posed in $H^1(\R^n)$. It is known that there exists solutions which form 
singularities in finite time. On the contrary, initial data with small $H^1(\R^n)$ norm 
yield global solutions. Furthermore, solitary waves are unstable. 

The effect of a noise on the behavior of the solutions has also been the object of several 
studies, both in the physical  literature (see for instance \cite{abdu-garnier}, 
\cite{BCIR1}, \cite{BCIR2}, 
\cite{FKLT}, \cite{konotop-vasquez}, \cite{UK}) or in the mathematical  literature
(see for instance  \cite{dBDL2}, \cite{dBDH1}, \cite{dBDblowupadditif},
\cite{dBDblowupmultip}, \cite{DDM}, \cite{deb-gautier},\cite{gautier1}, \cite{gautier2}). Random effects may be taken into 
account at various places of the equation. A random forcing term or a random potential can
be added. Also random diffraction index result as a random coefficient before the nonlinear
term. Numerical and theoretical studies have shown that many interesting new
behaviors may appear. 

For instance, 
it has been shown that  when a random potential which is white in time is added to the equation it 
may  affect strongly the formation of singularities. 
If this random potential is smooth in space and the nonlinearity is supercritical,
any initial data yields a solutions which blows up in finite time with 
positive probability. If the noise is additive, this is also true for critical nonlinearity. On the contrary, 
numerical experiments tend to show that, if the noise acts as a potential and is rough in space,
the formation of singularities is prevented and the solution continue to propagate. The rigorous
justification of such statement seems to be completely out of reach at present. 

In this work, we consider a noisy dispersion. This is a natural model in dispersion
managed optical fibers \cite{abp}, \cite{agarwal1}, \cite{agarwal2}, \cite{Ga}, \cite{marty} (see also 
 \cite{ZGJT01} for a deterministic periodic dispersion). The nonlinear Schr\"odinger equation with random 
dispersion has also been studied mathematically. In \cite{marty}, the power law nonlinearity
is replaced by a smooth bounded function and it is shown that, in a certain scaling, the solutions 
to the nonlinear Schr\"odinger equation converge to the solutions of the nonlinear Schr\"odinger
equation with white noise dispersion. This result has been extended to the case of a subcritical
nonlinearity in \cite{dBD10}. One of the main improvement in \cite{dBD10} is the use of Strichartz
type estimates for white noise dispersion (see also \cite{dB-R} for the derivation of Strichartz estimates for a stochastic Nonlinear Schr\"odinger equation). 

Note that Strichartz type estimates are not immediate for a white noise dispersion. 
We have an explicit formula of the fundamental solution for the linear  equation as in the deterministic case:
$$
u(t)= \frac1{\left(4i\pi\left(\beta(t)-\beta(s)\right)\right)^{d/2}}\int_{\R^d}\exp\left(\ds i\frac{|x-y|^2}{4(\beta(t)-\beta(s))} \right) u_s(y) dy
$$
is the solution of the linear equation with white noise dispersion with initial data $u_s$ at
time $s$ (see Proposition \ref{p3.1}).

Nevertheless, it is not obvious whether the Strichartz type estimate holds or not unlike the 
deterministic case. We have two difficulties to prove the Strichartz type estimate. One difficulty is that the dispersion coefficient is highly degenerate. In fact, for $\epsilon > s \ge 0$, the set
$\{t\in (s,\epsilon)\; : \beta(t)-\beta(s)=0\}$ has the cardinality of the continuum (see, e.g. 
\cite{durrett}, Example 4.1 in Section 7.4). Roughly speaking, in our problem, the dispersion coefficient 
has so many zeros that we can not expect that pathwise Strichartz estimates hold. 
Another difficulty is that the duality argument (or $T T^*$ argument) does not work as well as in the deterministic case. "Duality" corresponds to solving the equation backwards. For stochastic equations, a backward equation has in general no solution unless the coefficient of the noise
is considered as an unknown, which is not desirable in our situation. 

In the present work, we show that in the one dimensional case
it is possible to improve the Strichartz estimates obtained in \cite{dBD10} and as a result prove that the nonlinear 
Schr\"odinger equation with critical nonlinearity and white noise dispersion is globally well 
posed in $L^2(\R)$ and $H^1(\R)$. This confirms the fact that such a random dispersion has 
a strong stabilizing effect on the equation: in the quintic one dimensional case considered, it 
prevents the formation of singularities and yields global well posedness.

\section{Preliminaries and main results}

We consider the following stochastic nonlinear Schr\"odinger (NLS) equation with quintic nonlinearity 
on the real line
\begin{equation}
\label{e1.0}
\left\{
\begin{array}{l}
i du +\Delta u \circ d\beta + |u|^{4}u \,dt =0,\; x\in\R,\; t>0,\\
u(0)=u_0, \; x\in\R.
\end{array}
\right.
\end{equation}
The unknown $u$ is a random process on a probability space $(\O,\F,\P)$ depending
on $t>0$ and $x\in \R$. The noise term is given by a brownian motion
$\beta$  associated to a stochastic basis $(\Omega, \F,\P,(\F_t)_{t\ge 0})$.  
The product $\circ$ is a 
Stratonovich product. Classically, we transform this Stratonovitch equation into an It\^o equation which is formally equivalent:
\begin{equation}
\label{e1.1}
\left\{
\begin{array}{l}
\ds i du +\frac i2 \Delta^2 u \, dt+\Delta u \, d\beta  + |u|^{4}u \,dt=0,\; x\in\R,\; t>0,\\
u(0)=u_0.
\end{array}
\right.
\end{equation}
It seems as if the principal part of (2.2) were the double Laplacian, which does not appear to be degenerate.
But this is not true.
Indeed, the explicit formula of the fundamental solution for the linear equation shows the high degeneracy of the principal part (see Proposition 3.1), as is already pointed out in Section 1.

We study this equation (2.2) in the framework of the  $L^2(\R)$ based Sobolev spaces denoted 
by $H^s(\R)$, $s\ge 0$. 
We also use 
the spaces $L^p(\R)$ to treat the nonlinear term thanks to the Strichartz estimates. 
Note that, in all the article, these are spaces of complex valued functions.

For time dependent functions on an interval $I\subset \R$ with values in a Banach space $K$,
we use the spaces: $L^r(I; K), \, r\ge 1$. Given a time dependent function $f$, we use two notations for its values at some time $t$ depending on the context. We either write $f(t)$ or $f_t$.

The norm of a Banach space $K$ is simply denoted by $\|\cdot\|_K$. When we consider 
random variables with values in a Banach space $K$, we use $L^p(\O;K), \, p\ge 1$.

For spaces of predictible time dependent processes, we use the subscript $\mathcal P$. 
For instance
$L^r_{\mathcal P}(\O;L^p(0,T;K))$ is the subspace of $L^r(\O;L^p(0,T;K))$ consisting of
predictible processes.

Our  main result is the following.

\begin{Theorem}
\label{t2.1}
Let $u_0\in L^2(\R)$ a.s. be $\F_0$-measurable, then there 
exists a unique solution $u$  to 
\eqref{e1.1} with paths a.s. in $L^5_{loc}(0,\infty;L^{10}(\R))$; moreover, 
$u$ has paths  in $C(\R^+; L^2(\R))$, a.s. and
$$
\|u(t)\|_{L^2(\R)} =\|u_0\|_{L^2(\R)}, \quad a.s.
$$

If in addition $u_0 \in H^1(\R)$, then $u$ has paths a.s. in $C(\R^+;H^1(\R))$.
\end{Theorem}

As in \cite{dBD10}, we use this result to  
 justify rigorously the convergence of 
the solution of the  following random equation 
\begin{equation}
\label{e0.1}
\left\{
\begin{array}{l}
i \ds \frac{du}{dt} +\frac1{\vep}m\left(\frac{t}{\vep^2}\right) \partial_{xx} u   + |u|^{4}u =0,\; x\in\R,\; t>0,\\
u(0)=u_0,\, x\in\R,
\end{array}
\right.
\end{equation}
to the solution of \eqref{e1.1} provided that the real valued centered stationary random process $m(t)$ is continous
and that 
for any $T>0$, the process
$t \mapsto \vep \int_0^{t/\vep^2} m(s)ds$
converges in distribution to a standard real valued Brownian motion in $C([0,T])$. This is a classical 
assumption and can be verified in many cases. 

To our knowledge, Strichartz estimates are not available for equation \eqref{e0.1}. 
Hence we cannot get 
solutions in $L^2(\R)$. Since the equation is set in space dimension $1$, a local existence
result can be easily proved in $H^1(\R)$. For fixed $\vep$, we do not expect to have 
global in time solutions, indeed with a quintic nonlinearity it is known that singularities appear for 
the deterministic nonlinear Schr\"odinger equation. In the following result, we prove that the lifetime of the solutions
converges to infinity  when $\vep$  goes to zero, and that solutions of \eqref{e0.1} 
converge in distribution
to the solutions of the white noise driven equation (\ref{e1.1}).

\begin{Theorem}
\label{t2.2}
Suppose that $m$ satisfies the above assumption. Then, for any $\vep>0$ and 
$u_0\in H^1(\R)$, there exists a unique solution $u_\vep$ of equation
\eqref{e0.1} with continuous paths in $H^1(\R)$ which is defined on a random 
interval $[0,\tau_\vep(u_0))$. Moreover, for any $T>0$
$$
\lim_{\vep \to 0} \P(\tau_\vep(u_0)\le T)= 0,
$$
and the process $u_\vep\un_{[\tau_\vep > T]}$ 
converges in distribution to the solution $u$ of \eqref{e1.1} in $C([0,T];H^1(\R))$.
\end{Theorem}
\begin{Remark}
Note that there is a slight improvement compared to the result obtained in 
\cite{dBD10} where the convergence was not proved in the $H^1(\R)$ topology. This result can 
be extended to initial data in $H^s(\R)$ for $s\in (1/2,1]$. In this case, the convergence holds 
in $C([0,T];H^s(\R))$.
\end{Remark}

\section{The linear equation and Strichartz type estimates}

The Strichartz estimates are crucial to study the deterministic equation. In \cite{dBD10}, these
have been generalized to a white noise dispersion. However, the result obtained there
was not strong enough to treat the nonlinearity of the present article. We now show that
in dimension $1$, it is possible to get a better result.

We consider  the following stochastic linear Schr\"odinger equation:
\begin{equation}
\label{e3.2}
\left\{
\begin{array}{l}
\ds i du +\frac i2 \Delta^2 u \, dt+\Delta u \, d\beta  =0,\; t\ge s,\\
u(s)=u_s.
\end{array}
\right.
\end{equation}
We have an explicit formula for the solutions of 
\eqref{e3.2}. We recall from \cite{dBD10}, \cite{marty} the following result:
\begin{Proposition}
\label{p3.1}
For any $s\le T$ and $u_s\in {\mathcal S}'(\R^n)$, there exists a unique solution 
of \eqref{e3.2} almost
surely in $C([s,T];{\mathcal S}'(\R^n))$
and adapted. Its Fourier transform in space is given by 
$$
\hat u(t,\xi)= e^{-i|\xi|^2(\beta(t)-\beta(s))}\hat u_s(\xi), \; t\ge s,\; \xi \in \R^d.
$$
Moreover, if $u_s\in H^\s(\R)$ for some $\s\in \R$, then $u(\cdot)\in C([0,T];H^\s(\R))$ a.s. and 
$\|u(t)\|_{H^\s}=\|u_s\|_{H^\s}$,
a.s. for $t\ge s$.

If $u_s\in L^1(\R)$, the solution $u$ of \eqref{e3.2} has the expression 
\begin{equation}
\label{e3.3}
u(t)=S(t,s)u_s:= \frac1{\left(4i\pi\left(\beta(t)-\beta(s)\right)\right)^{d/2}}\int_{\R^d}\exp\left(\ds i\frac{|x-y|^2}{4(\beta(t)-\beta(s))} \right) u_s(y) dy, \; t\in [s,T].
\end{equation}
\end{Proposition}

\bigskip

The idea is to obtain Strichartz estimate through smoothing effects of $S(t,s)$ as was done 
in the deterministic case in \cite{ozawa-tsutsumi}.

The first step is the following.

\begin{Proposition}
\label{p3.2}
Let  $f\in L^4_{\mathcal P}(\O;L^1(0,T;L^2(R)))$ then 
$t\mapsto D^{1/2}\left(\left|\int_0^t S(t,s)f(s)ds\right|^2\right)$ belongs to $L^2_{\mathcal P}(\O\times [0,T]
\times \R)$ and
$$
\E\int_0^T\left\| D^{1/2}\left(\left|\int_0^t S(t,s)f(s)ds\right|^2\right)\right\|_{L^2(\R)}^2 dt
\le 4\sqrt{2 \pi}\;T^{1/2} \E\left( \| f\|^4_{L^1(0,T;L^2(\R))}\right).
$$
\end{Proposition}

\noindent
{\bf Proof.} By density, it is sufficient to prove that the inequality is valid for sufficiently smooth $f$.
Set, for $\xi\in\R$,
$$
A(\xi)=\left|\F\left[\left|\int_0^tS(t,s)f(s)ds\right|^2\right](\xi)\right|^2.
$$
Then, by Plancherel identity,
$$
\E\int_0^T\left\| D^{1/2}\left(\left|\int_0^t S(t,s)f(s)ds\right|^2\right)\right\|_{L^2(\R)}^2 dt
=\E\int_0^T\int_\R |\xi| A(\xi)d\xi dt.
$$
We have, by Proposition \ref{p3.1} and easy computations,
$$
\F\left[\left|\int_0^tS(t,s)f(s)ds\right|^2\right](\xi)
=\int_\R\int_0^t\int_0^t e^{-i(\beta_t-\beta_{s_1})(\xi-\xi_1)^2+i(\beta_t-\beta_{s_2})\xi_1^2}
\hat f_{s_1}(\xi-\xi_1)\hat{\bar f}_{s_2}(\xi_1)ds_1ds_2d\xi_1.
$$
We deduce:
$$
\ba{lr}
A(\xi)&=\ds \int\hspace{-.3cm}\int_{\R^2}\int\hspace{-.3cm}\int\hspace{-.3cm}\int\hspace{-.3cm}\int_{[0,t]^4} 
e^{-i(\beta_t-\beta_{s_1})(\xi-\xi_1)^2+i(\beta_t-\beta_{s_2})\xi_1^2} 
e^{i(\beta_t-\beta_{s_3})(\xi-\xi_2)^2-i(\beta_t-\beta_{s_4})\xi_2^2}\\
\\
& \ds\times \hat f_{s_1}(\xi-\xi_1)\hat{\bar f}_{s_2}(\xi_1)\bar{\hat f}_{s_3}(\xi-\xi_2)\bar{\hat {\bar f}}_{s_4}(\xi_2)
ds_1ds_2ds_3ds_4d\xi_1d\xi_2.
\ea
$$
Let us split $\ds[0,t]^4=\bigcup_{i=1,\dots,4} R_i$ with 
$$
R_i=\{(s_1,s_2,s_3,s_4)\in [0,t]^4;\; s_i=\max\{s_1,s_2,s_3,s_4\}\}
$$ 
and split accordingly 
$$
A(\xi)=\sum_{i=1,\dots,4}I_i(\xi).
$$
We then write, using $(\xi-\xi_1)^2-\xi_1^2-(\xi-\xi_2)^2+\xi_2^2=2\xi(\xi_2-\xi_1)$,
$$
\ba{ll}
\ds \E\left( \int_\R |\xi |ÊI_1(\xi) d\xi \right) 
= &\ds \E\bigg( \int_{R_1}\int\hspace{-.3cm}\int\hspace{-.3cm}\int_{\R^3} |\xi| 
e^{-2i(\beta_t-\beta_{s_1})\xi(\xi_2-\xi_1)-i(\beta_{s_1}-\beta_{s_2})\xi_1^2+i(\beta_{s_1}-\beta_{s_3})(\xi-\xi_2)^2-i(\beta_{s_1}-\beta_{s_4})\xi_2^2}\\
\\
& \ds\times \hat f_{s_1}(\xi-\xi_1)\hat{\bar f}_{s_2}(\xi_1)\bar{\hat f}_{s_3}(\xi-\xi_2)\bar{\hat {\bar f}}_{s_4}(\xi_2)
ds_1ds_2ds_3ds_4d\xi_1d\xi_2d\xi\bigg).
\ea 
$$
Clearly $\ds e^{-2i(\beta_t-\beta_{s_1})\xi(\xi_2-\xi_1)}$ is independent to the other factors. Moreover:
$$
\E\left( e^{-2i(\beta_t-\beta_{s_1})\xi(\xi_2-\xi_1)}\right) = e^{-2(t-s_1)\xi^2(\xi_2-\xi_1)^2}.
$$
We deduce
$$
\ba{ll}
\ds \E\left( \int_\R |\xi |ÊI_1(\xi) d\xi \right) \le&\ds \E\bigg( \int_{R_1}\int\hspace{-.3cm}\int\hspace{-.3cm}\int_{\R^3} |\xi| e^{-2(t-s_1)\xi^2(\xi_2-\xi_1)^2} |\hat f_{s_1}(\xi-\xi_1)||\hat{\bar f}_{s_2}(\xi_1)|\\
 \\
& \times |{\hat f}_{s_3}(\xi-\xi_2)||\hat{\bar f}_{s_4}(\xi_2)|
ds_1ds_2ds_3ds_4d\xi_1d\xi_2d\xi\bigg).
\ea 
$$
Note that
$$
\ba{l}
\ds\int\hspace{-.3cm}\int\hspace{-.3cm}\int_{\R^3} |\xi| e^{-2(t-s_1)\xi^2(\xi_2-\xi_1)^2} |\hat f_{s_1}(\xi-\xi_1)||\hat{\bar f}_{s_2}(\xi_1)| |{\hat f}_{s_3}(\xi-\xi_2)||\hat{\bar f}_{s_4}(\xi_2)|
d\xi_1d\xi_2d\xi\\
\ds= \int_\R |\xi| \left( \int_\R |\hat f_{s_1}(\xi-\xi_1)||\hat{\bar f}_{s_2}(\xi_1)| \left(\int_\R 
e^{-2(t-s_1)\xi^2(\xi_2-\xi_1)^2}|{\hat f}_{s_3}(\xi-\xi_2)||\hat{\bar f}_{s_4}(\xi_2)| d\xi_2\right)
d\xi_1\right) d\xi.
\ea
$$
Since $\ds \int_\R e^{-2(t-s_1)\xi^2\eta^2}d\eta= \frac{\sqrt{\pi}}{|\xi| \bigl (2(t-s_1) \bigr )^{1/2}}$, we deduce by Young's and Schwarz's inequalities:
$$
\ba{l}
\ds\int\hspace{-.3cm}\int\hspace{-.3cm}\int_{\R^3} |\xi| e^{-2(t-s_1)\xi^2(\xi_2-\xi_1)^2} |\hat f_{s_1}(\xi-\xi_1)||\hat{\bar f}_{s_2}(\xi_1)| |{\hat f}_{s_3}(\xi-\xi_2)||\hat{\bar f}_{s_4}(\xi_2)|
d\xi_1d\xi_2d\xi\\
\\
\ds\le  \frac{\sqrt{\pi}}{\bigl (2(t-s_1) \bigr )^{1/2}}\int_\R  \left( \int_\R |\hat f_{s_1}(\xi-\xi_1)|^2|\hat{\bar f}_{s_2}(\xi_1)|^2d\xi_1\right)^{1/2}\left(\int_\R 
|{\hat f}_{s_3}(\xi-\xi_2)|^2|\hat{\bar f}_{s_4}(\xi_2)|^2 d\xi_2\right)^{1/2}
 d\xi\\
 \\
 \dsÊ\le  \frac{\sqrt{\pi}}{\bigl ( 2(t-s_1) \bigr )^{1/2}}  \left( \int_\R\int_\R |\hat f_{s_1}(\xi-\xi_1)|^2|\hat{\bar f}_{s_2}(\xi_1)|^2d\xi_1d\xi\right)^{1/2}\left(\int_\R\int_\R 
|{\hat f}_{s_3}(\xi-\xi_2)|^2|\hat{\bar f}_{s_4}(\xi_2)|^2 d\xi_2d\xi \right)^{1/2}\\
\\
\ds=\frac{\sqrt{\pi}}{\bigl ( 2(t-s_1) \bigr )^{1/2}}  \|f_{s_1}\|_{L^2(\R)} \|f_{s_2}\|_{L^2(\R)} \|f_{s_3}\|_{L^2(\R)} \|f_{s_4}\|_{L^2(\R)}.
\ea
$$
It follows
$$
 \E\left( \int_\R |\xi |ÊI_1(\xi) d\xi \right) \le \E\int_{R_1}\frac{\sqrt{\pi}}{\bigl ( 2(t-s_1) \bigr )^{1/2}}
 \|f_{s_1}\|_{L^2(\R)} \|f_{s_2}\|_{L^2(\R)} \|f_{s_3}\|_{L^2(\R)} \|f_{s_4}\|_{L^2(\R)}ds_1ds_2ds_3ds_4
$$
and 
$$
\E\int_0^T\int_\R |\xi| I_1(\xi)d\xi dt \le \sqrt{2 \pi}T^{1/2} \E\left(\left( \int_0^T \|f_s\|_{L^2(\R)} ds \right)^4\right).
$$
The three other terms are treated similarly and the result follows. \hfill$\Box$

\begin{Proposition}
\label{p3.6}
There exists a constant 
$\kappa >0$ such that for any 
$s\in \R,\; T\ge 0$ and 
$f\in L_{\mathcal P}^4 (\O;L^{1}(s,s+T;L^{2}(\R)))$, the mapping 
$t\mapsto \int_s^t S(t,\s)f(\s)d\s$ belongs to $L_{\mathcal P}^4 (\O;L^{5}(s,s+T;L^{10}(\R)))$ and
$$
\left\|\int_s^\cdot S(\cdot,\sigma)f(\sigma)d\sigma\right\|_{L^4 (\O;L^{5}(s,s+T;L^{10}(\R)))}
\le \kappa T^{1/10} \left\|f\right\|_{L^4 (\O;L^{1}(s,s+T;L^{2}(\R)))}
$$
\end{Proposition}

\begin{Remark}
\label{r3.7}
This result is very similar to the classical Strichartz estimates in the case of  dimension $1$
considered here. Indeed $(5,10)$ and $(\infty,2)$ are admissible pairs. However, 
it is more powerful. Indeed, we have the extra factor $T^{1/10}$. This is a major 
difference and allows us to construct solution for the quintic nonlinearity. Recall that 
in the deterministic case, it is known that there are singular solutions for this equation. 
The proof below extends easily to the same result with $(5,10)$ replaced by any 
admissible pair $(r,p)$, {\it i.e.} satisfying $\frac2r=\frac12-\frac1p$. Of course, the power
of $T$ changes in this case; but it remains positive.
\end{Remark}

{\bf Proof:} We treat only the case $s=0$. The generalization is easy. Also, it is sufficient to prove that the
inequality holds for sufficiently smooth $f$. 

We use the following Lemma. Its proof is given below for the reader's convenience.
\begin{Lemma}
\label{l3.6}
Let $g\in L^1(\R)$ such that $D^{1/2}g\in L^2(\R)$, then $g\in L^5(\R)$ and
$$
\|g\|_{L^5(\R)}\le C\|g\|_{L^1(\R)}^{1/5}\|D^{1/2}g\|_{L^2(\R)}^{4/5}.
$$
\end{Lemma}

Let us write
$$
\left\|\int_0^\cdot S(\cdot,\sigma)f(\sigma)d\sigma\right\|_{L^4 (\O;L^{5}(0,T;L^{10}(\R)))}^4=
\left\|\left|\int_0^\cdot S(\cdot,\sigma)f(\sigma)d\sigma\right|^2\right\|_{L^2 (\O;L^{5/2}(0,T;L^{5}(\R)))}^2.
$$
Therefore, by Lemma \ref{l3.6}, H\"older inequality and Proposition \ref{p3.2}, 
$$
\ba{l}
\ds\left\|\int_0^\cdot S(\cdot,\sigma)f(\sigma)d\sigma\right\|_{L^4 (\O;L^{5}(0,T;L^{10}(\R)))}^4\\
\\
\ds\le c\E\left(\left( \int_0^T \left\|\left|\int_0^tS(t;\sigma)f_\sigma d\sigma\right|^2\right\|^{1/2}_{L^1(\R)}
 \left\|D^{1/2}\left|\int_0^tS(t;\sigma)f_\sigma d\sigma\right|^2\right\|_{L^2(\R)}^{2} dt \right)^{4/5}\right)\\
\\
\ds\le c\E\left(  \left\|\left|\int_0^\cdot S(\cdot;\sigma)f_\sigma d\sigma\right|^2\right\|^{2/5}_{L^\infty(0,T;L^1(\R))}
 \left\|D^{1/2}\left|\int_0^\cdot S(\cdot;\sigma)f_\sigma d\sigma\right|^2\right\|_{L^2(0,T;L^2(\R))}^{8/5}  \right)\\
 \\
\ds  \le c \E\left(  \left\|\left|\int_0^\cdot S(\cdot;\sigma)f_\sigma d\sigma\right|^2\right\|^{2}_{L^\infty(0,T;L^1(\R))}\right)^{1/5}
\E\left( \left\|D^{1/2}\left|\int_0^\cdot S(\cdot;\sigma)f_\sigma d\sigma\right|^2\right\|_{L^2(0,T;L^2(\R))}^{2}  \right)^{4/5}\\
\\
\ds \le T^{2/5} \E\left( \|f\|_{L^1(0,T;L^2(\R))}^4\right).
 \ea
$$
\hfill$\Box$

{\bf Proof of Lemma \ref{l3.6}:} By Gagliardo-Nirenberg inequality, we have:
\be
\label{e3.4}
\|g\|_{L^5(\R)}\le c \|D^{1/2}g\|_{L^2(\R)}^{3/5}\|g\|_{L^2(\R)}^{2/5}.
\ee
Moreover
$$
\ba{ll}
\ds \|g\|_{L^2(\R)}^2=\|\hat g\|_{L^2(\R)}^2&\ds =\int_{|\xi| \ge R} |\hat g(\xi)|^2 d\xi +\int_{|\xi| \le R} |\hat g(\xi)|^2 d\xi \\
\\
&\ds \le \int_{|\xi| \ge R} \frac{|\xi|}{R} |\hat g(\xi)|^2 d\xi + 2R \|\hat g\|_{L^\infty(\R)}^2\\
\\
&\ds\le \frac1R \|D^{1/2}g\|_{L^2(\R)}^2 + 2R \|g\|_{L^1(\R)}^2
\ea
$$
It suffices to take $R= \|D^{1/2}g\|_{L^2(\R)} \|g\|_{L^1(\R)}^{-1}$ and to insert the result in 
\eqref{e3.4} to conclude. \hfill $\Box$

We also need to have estimates on the action of $S(t,s)$ on an initial data.
\begin{Proposition}
\label{p3.6bis}
Let $s\ge 0$ and $u_s\in L^4(\O;L^2(\R))$ be $\F_s$ measurable, then $t\mapsto  S(t,s)u_s$ 
belongs to $L^4_{\mathcal P}(\Omega;L^5(s,s+T;L^{10}(\R)))$ and 
$$
\|S(\cdot,s)u_s\|_{L^4(\Omega;L^5(s,s+T;L^{10}(\R)))}\le c T^{1/10} \|u_s\|_{L^4(\O;L^2(\R))}.
$$
\end{Proposition}
{\bf Proof:} The proof is similar. Again, we only treat the case $s=0$. We first write:
$$
\left|\F\left(\left|S(t,0)u_0\right|^2\right)\right|^2=\int\hspace{-.3cm}\int_{\R^2}
e^{-2i\beta_t\xi(\xi_2-\xi_1)}{\hat u_0}(\xi-\xi_1)\hat{\bar u}_0(\xi_1)\bar{\hat u}_0(\xi-\xi_2)
\bar{\hat {\bar u}}_0(\xi_2)d\xi_1d\xi_2
$$
and
$$
\ba{l}
\ds \E\left( \left\| D^{1/2}\left|S(t,0)u_0\right|^2\right\|_{L^2(0,T;L^2(\R))}^2\right)\\
\ds=\E \int_0^T\int\hspace{-.3cm}\int\hspace{-.3cm}\int_{\R^3}|\xi|e^{-2t\xi^2(\xi_2-\xi_1)^2}
{\hat u_0}(\xi-\xi_1)\hat{\bar u}_0(\xi_1)\bar{\hat u}_0(\xi-\xi_2)
\bar{\hat {\bar u}}(\xi_2)d\xi_1d\xi_2d\xi dt\\
\\
\ds\le \E\int_0^T \int_\R |\xi| \left(\int_\R {\hat u_0}(\xi-\xi_1)\hat{\bar u}_0(\xi_1)
\left( \int_\R e^{-2t\xi^2(\xi_2-\xi_1)^2} \bar{\hat u}_0(\xi-\xi_2)
\bar{\hat {\bar u}}(\xi_2)d\xi_2\right)d\xi_1\right)d\xi dt.
\ea
$$
Therefore by Young's and Schwarz's inequalities:
$$
\ba{ll}
\ds \E\left( \left\| D^{1/2}\left|S(t,0)u_0\right|^2\right\|_{L^2(0,T;L^2(\R))}^2\right)&
\ds\le \E\int_0^T  \sqrt{\pi} t^{-1/2} \E\left(\|u_0\|_{L^2(\R)}^4\right) dt\\
&\ds \le 2\sqrt{\pi}T^{1/2} \E\left(\|u_0\|_{L^2(\R)}^4\right) .
\ea
$$
We then use Lemma \ref{l3.6} and H\"older inequality:
$$
\ba{l}
\ds \|S(\cdot,0)u_0\|_{L^4(\Omega;L^5(0,T;L^{10}(\R)))}\\
\\
\ds \le c
\left\|\left|S(\cdot,0,u_0)\right|^2\right\|_{L^2(\O;L^\infty(0,T;L^1(\R)))}^{1/10}
\left\|D^{1/2}\left|S(\cdot,0,u_0)\right|^2\right\|_{L^2(\O;L^2(0,T;L^2(\R)))}^{4/10}\\
\\
\dsÊ\le cT^{1/10} \E\left(\|u_0\|_{L^2(\R)}^4\right) .
\ea
$$
\hfill $\Box$
\section{Proof of Theorem \ref{t2.1}}

As is classical, we first construct a local solution of equation \eqref{e1.1} thanks to a cut-off of the nonlinearity. Proceeding as in 
in \cite{dBDL2}, \cite{dBDH1}, \cite{dBD10}, we take
 $\theta \in C_0^\infty(\R)$ be such that $\theta =1$ on $[0,1]$, 
$\theta=0$ on $[2,\infty )$ and for $s\in \R$,
$u\in L^5_{loc}(s,\infty;L^{10}(\R))$, $R\ge 1$ and $t\ge 0$, we set
$$
\theta^s_R(u)(t)=\theta\left( \frac{\|u\|_{L^5(s,s+t;L^{10}(\R))} }R\right).
$$
For $s=0$, we set $\theta^0_R=\theta_R$.

The truncated form of equation \eqref{e1.1} is given by
\begin{equation}
\label{e4.2}
\left\{
\begin{array}{l}
\ds i du^R +\frac{i}2 \Delta^2u^R\, dt+\Delta u^R  d\beta +\theta_R(u^R) |u^R|^{4}u^R \,dt =0,\\
u^R(0)=u_0.
\end{array}
\right.
\end{equation}
We interpret it in the mild sense
\be
\label{e4.3}
u^R(t)=S(t,0)u_0 +i \int_0^t S(t,s)\theta_R(u^R)(s) |u^R(s)|^{4} u^R(s) ds.
\ee

\begin{Proposition}
\label{t4.3}
For any $\F_0$-measurable $u_0\in L^4(\O;L^2(\R))$, 
there exists a unique solution of 
\eqref{e4.3} $u^R$ in $L^4_{\mathcal P}(\O;L^5(0,T;L^{10}(\R))))$ for any $T>0$. Moreover $u^R $  is a weak solution of \eqref{e4.2} in the sense that
for any $\varphi\in C_0^\infty(\R^d)$ and any $t\ge 0$,
$$
\begin{array}{rl}
& i(u^R(t)-u_0,\varphi)_{L^2(\R)}
\\
= & \ds -\frac{i}2\int_0^t (u^R,\Delta^2\varphi)_{L^2(\R)}ds - \int_0^t
\theta_R(u^R)(|u^R|^{4}
u^R,\varphi)_{L^2(\R)}ds - \int_0^t (u^R,\Delta \varphi)_{L^2(\R)}d\beta(s),\; a.s.
\ea
$$
Finally, the $L^2(\R)$ norm is conserved:
$$
\|u^R(t)\|_{L^2(\R)}=\|u_0\|_{L^2(\R)},\; t\ge 0, \; a.s.
$$
and $u\in 
C([0,T];L^2(\R))$ a.s.
\end{Proposition}

\noindent
{\bf Proof.} In order to lighten the  notations we omit the $R$ dependence in this 
proof. By Proposition \ref{p3.6bis}, we know that 
$S(\cdot,0)u_0\in L^4_{\mathcal P}(\O;L^5(0,T;L^{10}(\R))))$. Then, by Proposition \ref{p3.6}, for $u,v\in L^4_{\mathcal P}(\O;L^5(0,T;L^{10}(\R))))$,
$$
\ba{l}
 \ds\left\| \int_0^t S(t,s)\left( \theta(u)(s) |u(s)|^{4} u(s)
-\theta(v)(s) |v(s)|^{4} v(s)\right) ds \right\|_{L^4(\O;L^5(0,T;L^{10}(\R))))}\\
\\
 \ds \le c T^{1/10} \left\| \theta(u) |u|^{4} u-\theta(v) |v|^{4} v \right\|_{L^4(\O;L^1(0,T;L^2(\R)))}\\
 \\
\le c\, T^{1/10} R^{4}\|u-v\|_{L^4(\O;L^5(0,T;L^{10}(\R))))}.
\ea
$$
It follows that 
\be
\label{cont}
\T^R : u\mapsto S(t,0)u_0 +i \int_0^t S(t,s)\theta(u(s)) |u(s)|^{4} u(s) ds
\ee
defines a strict contraction on $L^4_{\mathcal P}(\O;L^5(0,T;L^{10}(\R))))$ provided 
$T\le T_0$ where $T_0$ depends only on $R$. Iterating this construction, one easily 
ends the proof of the first statement. 
The proof that $u$ is in fact a weak solution is classical. 

Let $M\ge 0$ and $u_M=P_M u $ be a regularization of the solution $u$ defined by a 
truncation in Fourier space: $\hat u_M(t,\xi)=\theta\left(\frac{|\xi|}{M}\right) \hat u(t,\xi)$. We deduce from the weak form of the equation that
$$
idu_M +\frac{i}2 \Delta^2u_M\, dt+\Delta u_M  d\beta +P_M \left(\theta(u) |u|^{4}u\right) \,dt =0.
$$
We apply It\^o formula to $\|u_M\|^2_{L^2(\R)}$ and obtain
$$
\|u_M(t)\|_{L^2(\R)}^2=\|u_0\|^2_{L^2(\R)}+ Re\left( i\int_0^t  \left(\theta(u) |u|^{4}u, u_M\right)ds
\right),\; t\in [0,T].
$$
We know that $u\in L^5(0,T;L^{10}(\R))$ a.s. By the integral equation,
$$
\ba{ll}
\|u(t)\|_{L^2(\R)}&\ds \le \|S(t,0)u_0\|_{L^2(\R)} + \int_0^t \|S(t,s)\theta(u(s)) |u(s)|^{4} u(s) \|_{L^2(\R)}ds\\
&\ds \le \|u_0\|_{L^2(\R)} + \int_0^t \| u(s) \|_{L^{10}(\R)}^5ds.
\ea
$$
We deduce that $u\in L^\infty(0,T;L^2(\R))$ a.s. and 
$$
\lim_{M\to \infty}  u_M =u \mbox{ in } L^\infty(0,T;L^2(\R)),\; \mbox{a.s.}
$$
we may let $M$ go to infinity in the above equality and obtain 
$$
\lim_{M\to \infty} \|u_M(t)\|_{L^2(\R)}=\|u_0\|_{L^2(\R)},\; t\in [0,T],\; \mbox{a.s.}
$$
This implies $u(t)\in L^2(\R)$ for  any $t\in [0,T]$ and $\|u(t)\|_{L^2(\R)}=\|u_0\|_{L^2(\R)}$. 
As easily seen 
from the weak form of the equation, $u $ is almost surely continuous with values in 
$H^{-4}(\R)$. It follows that $u$ is weakly continuous with values in $L^2(\R)$. 
Finally the continuity of $t\mapsto \|u(t)\|_{L^2(\R)}$ implies $u\in 
C([0,T];L^2(\R))$. 
\hfill $\square$

\bigskip

The construction of a global solution and the end of the proof of Theorem \ref{t2.1} are
 now very similar to what was done in \cite{dBD10}.
We briefly recall the ideas for the reader's convenience.

There is no loss of generality in assuming that $u_0\in L^2(\R)$ is deterministic.
Uniqueness is clear since two 
solutions are solutions of the truncated equation on a random interval. We fix $T_0$ and construct a 
solution on $[0,T_0]$.

We define 
$$
\tau_R=\inf\{ t\in [0,T],\; \|u^R\|_{L^5(0,t; L^{10}(\R))}\ge R  \}
$$
so that $u^R$ is a solution of \eqref{e1.1} on $[0,\tau_R]$. 
\begin{Lemma}
\label{l5.1}
There exist constants $c_1,\, c_2$ such that if 
$$
T^{2/5} \le c_1\; R^{-16} 
$$
then
$$
\P(\tau_R \le T) \le \frac{c_2\|u_0\|^4_{L^2(\R)}}{R^4}
$$
\end{Lemma}

\noindent
{\bf Proof.} We write
\be
\label{eqint}
u^R(t)\un_{[0,\tau_R]}(t)= S(t,0)u_0\un_{[0,\tau_R]}(t) + i 
\int_0^t S(t,s)|u^R|^{4}u^R \un_{[0,\tau_R]}(s)ds \un_{[0,\tau_R]}(t).
\ee
Thus for $T\le T_0$ 
$$
\ba{ll}
\ds \|u^R\un_{[0,\tau_R]}\|_{L^5(0,T;L^{10}(\R))}& \ds 
\le \| S(\cdot,0)u_0\un_{[0,\tau_R]}\|_{L^5(0,T;L^{10}(\R))} \\
&\ds + \|\int_0^t S(t,s)|u^R|^{4}u^R \un_{[0,\tau_R]}(s)ds \|_{L^5(0,T;L^{10}(\R))}.
\ea
$$
Proposition \ref{p3.6} and Proposition \ref{p3.6bis} yield
$$
\ba{ll}
\ds \E\left(\|u^R\un_{[0,\tau_R]}\|_{L^5(0,T;L^{10}(\R))}^4\right)&\ds \le c(T_0) \|u_0\|_{L^2(\R)}^4+
c \, T^{2/5} \E\left(\left\|\left|u^R\right|^{5}\un_{[0,\tau_R)}\right\|_{L^1(0,T;L^{2}(\R))}^4\right)\\
\\
&\ds  \le  \ds c(T_0) \|u_0\|_{L^2(\R)}^4+
c \, T^{2/5} \E\left(\left\|u^R\un_{[0,\tau_R)}\right\|_{L^5(0,T;L^{10}(\R))}^{20}\right)\\
\\
&\ds  \le  \ds c(T_0) \|u_0\|_{L^2(\R)}^4+
c \, T^{2/5} R^{16}\E\left(\left\|u^R\un_{[0,\tau_R)}\right\|_{L^5(0,T;L^{10}(\R))}^{4}\right)
\ea 
$$
Hence, if $c \,  T^{2/5} R^{16} \le \ds\frac12$,
$$
\E\left(\|u^R\un_{[0,\tau_R]}\|_{L^5(0,T;L^{10}(\R))}^4\right)
\ds \le 2 c(T_0) \|u_0\|_{L^2(\R)}^4
$$
and by Markov inequality
$$
\P(\tau_R\le T) \le \frac{2c(T_0)\|u_0\|^4_{L^2(\R)}}{R^4}.
$$
\hfill $\square$

\bigskip

In order to construct a solution to \eqref{e1.1} on $[0,T_0]$, we iterate the local construction. We fix 
$R>0$ and have a local solution on $[0,\tau_R]$.  We set $\tau^0_R=\tau_R$.
We then consider recursively the equation for $u$. For $n\ge 0$, we set $T_R^n=\ds\sum_{k=0}^n
\tau^n_R$ and define :
$$
u(t+T^n_R)=S(t+T^n_R,T^n_R)u(T^n_R) +\int_0^t S(t+T^n_R,s+T^n_R)\theta_R^{T^n_R}(u)(s) |u(s+T^n_R)|^{2\s} u(s+T^n_R) ds.
$$ 
The local construction can be reproduced and we obtain a unique global 
solution of this equation on $[T^n_R,T^n_R+\tau^{n+1}_R]$ where
$$
\tau_R^{n+1}=\inf\{ t\in [0,T],\; |u|_{L^5(T^n_R,t+T^n_R; L^{10}(\R))}\ge R  \}.
$$
We thus obtain a solution of the non truncated equation on $\ds\left[0,\sum_{n=0}^\infty \tau_R^n\right]$.
By Lemma \ref{l5.1}, the strong Markov property and the conservation of the $L^2(\R)$ norm 
$$
\P(\tau_R^{n+1}\le T|\F_{T^n_R})=\P(\tau_R^{n+1}\le T|u(T^n_R)) \le  \frac{c_2|u(T^{n}_R)|^4_{L^2(\R)}}{R^4} = \frac{c_2|u_0|^4_{L^2(\R)}}{R^4},\; a.s.,
$$
provided  $T^{2/5} \le c_1\; R^{-16} $.
Note that 
$$
\P\left(\lim_{n \to +\infty}\tau_R^n =0\right) = \lim_{\vep \to 0} \lim_{N\to +\infty} \P( \tau_R^n\le \vep,\;  n\ge N).
$$
Finally we choose $R$ large enough and $\vep^{2/5} \le c_1\; R^{-16} $ so that, for 
all $n\in\N$,
$$
\P( \tau_R^{n+1}\le \vep|\F_{T^{n}_R})\le \frac12,\; a.s.
$$
Then, since $\P( \tau_R^{M}\le \vep|\F_{T^{M-1}_R})=\E\left(\un_{\tau_R^M\le \vep}|\F_{T^{M-1}_R}\right)$, we have for  $0\le N\le M$:
$$
\ba{ll}
\ds \P( \tau_R^n \le \vep,\; M\ge n\ge N)&\ds =\E\left(\prod_{M\ge n\ge N}\un_{\tau_R^n\le \vep}\right)\\
&\ds 
=\E\left(\E\left(\un_{\tau_R^M\le \vep}|\F_{T^{M-1}_R}\right)\prod_{M-1\ge n\ge N}\un_{\tau_R^n\le \vep}\right)\\
\\
&\ds \le \frac12 \E\left(\prod_{M-1\ge n\ge N}\un_{\tau_R^n\le \vep}\right).
\ea
$$
Repeating the last inequality, we deduce
$$
\P( \tau_R^n \le \vep,\; M\ge n\ge N)\le \frac1{2^{M-N}}
$$
and
$$
\P( \tau_R^n \le \vep,\;  n\ge N)\le  \lim_{M\to \infty} \P( \tau_R^n \le \vep,\; M\ge n\ge N) \le \lim_{M\to \infty} \frac{1}{2^{M-N}}=0.
$$
Hence, $\P(\lim_{n\to +\infty}\tau_R^n =0) =0$ so that $\tau^0_R+\dots+\tau^n_R$ goes to infinity
a.s. and we have constructed a global solution.

The conservation of the $L^2$-norm and the fact that $u \in C(\R^+;L^2(\R))$ a.s. was proved in Theorem \ref{t4.3}.

Finally, assume that $u_0 \in H^1(\R)$. Then going back to $\T^R$ defined in \eqref{cont}, and applying the same estimates
as in the proof of Lemma \ref{l5.1}, after having taken first order space derivatives, lead to
$$
\ba{l}
\|\T^Ru\|_{L^4(\O;L^5(0,T;W^{1,10}(\R))}\\ \\
\le CT_0^{1/10}\|u_0\|_{H^1(\R)}+C'T^{1/10}R^{16}\|u\|_{L^4(\O;L^5(0,T;W^{1,10}(\R))}
\ea
$$
This proves that   $B=B(0,R_0)$, the ball of radius $R_0$ in $L^4(\O;L^5(0,T;W^{1,10}(\R))$ is 
invariant by $\T^R$ provided $T\le \tilde T_0$, where $\tilde T_0$ depends only on $R$ and not on $R_0$.
Since closed balls of $L^4(\O;L^5(0,T;W^{1,10}(\R))$ are closed in $L^4(\O;L^5(0,T;W^{1,10}(\R))$, this implies 
that the fixed point of $\T^R$, which is the solution $u^R$ of \eqref{e4.3}, is in $L^4(\O;L^5(0,T;W^{1,10}(\R))$.

We deduce that $u$ has paths in $L^5(0,T_0;W^{1,10}(\R)$ and $|u|^4u$ in $L^1(0,T_0;H^1(\R))$.

It is easily proved that $t\mapsto \int_0^t S(t,s)f(s) ds$ is in 
$L^p(\O;C([0,T];H^1(\R))$ provided $f\in L^p(\O;L^1(0,T;H^1(\R)))$ and that 
$t\mapsto S(t,0)u_0$ is in $L^p(\O;C([0,T];H^1(\R))$ for $u_0\in L^p(\O;H^1(\R))$.

By a localization argument, we conclude that $u$ is continuous with values in $H^1(R)$ for
$u_0\in H^1(\R)$
\hfill $\square$

\section{Equation \eqref{e1.0} as limit of NLS equation with random dispersion}

The proof of Theorem \ref{t2.2} uses similar arguments as in \cite{dBD10}, however there are some modifications which enable us to get a stronger result. We fix $T\ge 0$.

Consider the following nonlinear Schr\"odinger equation written in the mild form:
$$
u_n(t)=S_n(t)u_0 +i \int_0^t S_n(t,\s) F(|u(\s)|^2)u(\s) d\s,
$$
where $F$ is a smooth function with compact support,
$n$ is a real valued function
and we have denoted by $S_n(t,\s)=\F^{-1} e^{-i(n(t)-n(\s))\xi^2/2}\F$,
the evolution operator associated to the linear 
equation
$$
i \ds \frac{dv}{dt} +\dot n(t) \partial_{xx} v    =0, \; x\in\R,\; t>0.
$$
Since $S_n(t,\s)$
is an isometry on  $H^1(\R)$, it is easily shown that for $u_0\in H^1(\R)$ there exists a unique $u_n$ 
in $C([0,T];H^1(\R))$, provided that
$n$ is a continuous function of $t$.

Let $(n_k)$ be a sequence in $C([0,T];\R)$ which converges to $n\in C([0,T];\R)$ uniformly on $[0,T]$.
Then, for $u_0\in H^1(\R)$, we have
$$
\ba{ll}
\ds \|u_{n_k}(t)-u_n(t)\|_{H^1(\R)}&\ds  \le \left\|\left(S_{n_k}(t,0)-S_n(t,0)\right)u_0\right\|_{H^1(\R)}
\\
\\
&\ds + \int_0^t \left\|\left(S_{n_k}(t,\s)-S_n(t,\s)\right)F(|u_n(\s)|^2)u_n(\s)\right\|_{H^1(\R)}d\s\\
\\
&\ds + \int_0^t \left\|S_{n_k}(t,\s)\left(F(|u_n(\s)|^2)u_n(\s)-F(|u_{n_k}(\s)|^2)u_{n_k}(\s)\right)\right\|_{H^1(\R)}d\s
\ea
$$ 
Since $F$ is smooth and has compact support, there exists $M_F$ such that
$$
\ba{ll}
\ds \|F(|u|^2)u-F(|v|^2)v\|_{H^1(\R)} & \ds \le M_F \left(\|u-v\|_{H^1(\R)} +\|u\|_{H^1(\R)}\| u-v\|_{L^\infty(\R)}\right)\\
& \ds \le M_F \left(\|u-v\|_{H^1(\R)} +\|u\|_{H^1(\R)}\| u-v\|_{H^1(\R)}\right).
\ea
$$
Since $S_{n_k}(t,\s)$ is an isometry, we deduce
$$
\ba{c}
\ds \int_0^t \left\|S_{n_k}(t,\s)\left(F(|u_n(\s)|^2)u_n(\s)-F(|u_{n_k}(\s)|^2)u_{n_k}(\s)\right)\right\|_{H^1(\R)}d\s \\
\ds \le C \int_0^t \left\|u_n(\s)-u_{n_k}(\s)\right\|_{H^1(\R)}d\s
\ea 
$$
with $C=M_F\left(1+\sup_{t\in[0,T]}\|u_n(t)\|_{H^1(\R)}\right)$.
It is easily checked that 
\be
\label{e5.1}
\left\|\left(S_{n_k}(t,0)-S_n(t,0)\right)u_0\right\|_{H^1(\R)}\to 0
\ee
as $k\to \infty$. Finally, note that $\{u_n(\s);\; \s\in [0,T]\}$ is compact in $H^1(\R)$. By continuity 
of $u\mapsto F(|u|^2)u$ on $H^1(\R)$, we deduce that $\{F(|u_n(\s)|^2)u_n(\s);\; \s\in [0,T]\}$ is 
also compact in $H^1(\R)$. It follows that for any $\delta$, we can find an $R_\delta$ such that
$$
\sup_{\s\in[0,T]} \left\| |\xi| \F\left( F(|u_n(\s)|^2)u_n(\s)\right)1_{|\xi|\ge R_\delta}\right\| _{L^2(\R)}\le \delta.
$$
Moreover, there exists $N_\delta\in\N$ such that, for $k\ge N_\delta$,
$$
\sup_{0\le \s\le t\le T} \left\||\xi|\left( e^{-i(n(t)-n(s))\xi^2/2}-e^{-i(n_k(t)-n_k(s))\xi^2/2}\right)
\F\left(F(|u_n(\s)|^2)u_n(\s)\right)1_{|\xi|\le R_\delta}\right\|_{L^2(\R)} \le \delta.
$$
We deduce
$$
 \int_0^t \left\|\left(S_{n_k}(t,\s)-S_n(t,\s)\right)F(|u_n(\s)|^2)u_n(\s)\right\|_{H^1(\R)}d\s
 \le 3T\delta
$$
for $k\ge N_\delta$. By \eqref{e5.1}, we may assume that 
$$
\left\|\left(S_{n_k}(t,0)-S_n(t,0)\right)u_0\right\|_{H^1(\R)}\le \delta
$$
for $k\ge N_\delta$. By Gronwall Lemma, we finally prove
$$
\sup_{t\in [0,T]}\|u_{n_k}(t)-u_n(t)\|_{H^1(\R)} \le (3T+1)e^{CT} \delta.
$$
This proves that the map $n\to u_n$ is continuous form $C([0,T])$ into 
$C([0,T];H^1(\R))$. 

Under our assumption, the process 
$t\mapsto \int_0^t \frac1\vep m(\frac{s}{\vep^2})ds$ converges in distribution in 
$C([0,T])$ to a brownian motion, and so we deduce that the solution of 
\begin{equation}
\label{e6.2}
\left\{
\begin{array}{l}
i \ds \frac{du}{dt} +\frac1\vep m(\frac{t}{\vep^2}) \partial_{xx} u   
+ F(|u|^{2})u =0, \; x\in\R,\; t>0,\\
\\
u(0)=u_0, \; x\in\R,
\end{array}
\right.
\end{equation}
converges in distribution in $C([0,T];H^1(\R))$ to the solution of 
$$
\left\{
\begin{array}{l}
i du +\Delta u \circ d\beta +F( |u|^{2})u \,dt =0,\; x\in\R,\; t>0,\\
u(0)=u_0,\; x\in\R.
\end{array}
\right.
$$
We now want to extend this result to the original power nonlinear term. 
Let us introduce the truncated equations, where $\theta$ is as in section 4,
\begin{equation}
\label{e6.3}
\left\{
\begin{array}{l}
\ds i \ds \frac{du}{dt} +\frac1\vep m(\frac{t}{\vep^2}) \partial_{xx} u   
+ \theta\left(\frac{|u|^2}{M}\right)|u|^{4}u =0, \; x\in\R,\; t>0,\\
u(0)=u_0, \;x\in\R,
\end{array}
\right.
\end{equation}
and
\begin{equation}
\label{e6.4}
\left\{
\begin{array}{l}
\ds i du +\Delta u \circ d\beta + \theta\left(\frac{|u|^2}{M}\right)|u|^{4}u \,dt =0,\; x\in\R,\; t>0,\\
u(0)=u_0,\; x\in\R.
\end{array}
\right.
\end{equation}
We denote by $u^M_\vep$ and $u^M$ their respective solutions. By the previous
arguments, these solutions exist and
are unique in $C([0,T];H^1(\R))$. Note that setting
$$
\widetilde\tau_\vep^M=\inf\{t\ge 0:\; \|u_\vep^M(t)\|_{L^\infty(\R)}\ge M\}
$$
and $u_\vep=u_\vep^M$  
on $[0,\widetilde\tau_\vep^M]$, defines a unique local solution $u_\vep$
of equation \eqref{e0.1} on $[0,\tau_\vep)$ with $\tau_\vep =\lim_{M\to \infty}
\tilde\tau_\vep^M$.

We also set
$$
\widetilde\tau^M=\inf\{t\ge 0:\; \|u^M(t)\|_{L^\infty(\R)}\ge M\}.
$$
By the above result, for each $M$, $u^M_\vep$
converges to  $u^M$ in distribution in $C([0,T];H^1(\R))$. 
By Skorohod 
Theorem, after a change of probability space, we can assume that for each $M$ 
the convergence of $u^M_\vep$ to $u^M$ holds almost surely in $C([0,T];H^1(\R))$.
To conclude, let us notice that for $0<\delta\le 1$, if 
$$
\widetilde\tau^{M-1}\ge T \mbox{ and } \|u_\vep^M-u^M\|_{C([0,T];H^1(\R))}\le \delta
$$
then $u^M=u$, the solution of \eqref{e1.1}, on $[0,T]$. Moreover, by the Sobolev 
embedding
$H^1(\R)\subset L^\infty(\R)$, we have 
$$
\|u_\vep^M-u^M\|_{C([0,T];L^\infty(\R))}\le c\delta
$$
for some $c>0$. We deduce $|u_\vep^M|_{C([0,T];L^\infty(\R))}\le M$ provided $\delta$ is small 
enough.
Therefore 
$$
\tau_\vep>\widetilde\tau_\vep^M\ge T \mbox{ and } u_\vep^M=u_\vep \mbox{ on }
[0,T].
$$ 
It follows that for $\delta>0$ small enough,
$$
\ba{rl}
& \P(\tau_\vep(u_0)\le T) +
\P(\tau_\vep(u_0)> T \mbox{ and } \|u_\vep-u\|_{C([0,T];H^1(\R))}>\delta ) \\
\\
\le & \P(\|u_\vep^M-u^M\|_{C([0,T];H^1(\R))}>\delta)
+\P(\tilde\tau^{M-1}<T).
\ea
$$
Since $u_0\in H^1(\R)$, we know that $u$ is almost surely in $C(\R^+;H^1(\R))$;
we deduce
$$
\lim_{M\to \infty} \P(\tilde\tau^{M-1}<T) =0.
$$
Choosing first $M$ large and then $\vep$ small we obtain
$$
\lim_{\vep\to 0} \P(\tau_\vep(u_0)\le T) = 0
$$
and
$$
\lim_{\vep \to 0}
\P(\tau_\vep(u_0)> T \mbox{ and } \|u_\vep-u\|_{C([0,T];H^1(\R))}>\delta )= 0
$$
The result follows.
\hfill $\square$


\begin{thebibliography}{99}

\bibitem{abp}F.Kh. Abdullaev, J.C. Bronski, G. Papanicolaou,
{\em Soliton perturbations
and the random Kepler problem},
{Physica D}
{\bf 135}, 369--386  (2000).

\bibitem{abdu-garnier} F. Kh. Abdullaev, J. Garnier,  
{\em Optical solitons in random media}
Progress in Optics, Vol. 48, pp. 35--106 (2005). 

\bibitem{agarwal1} G.P. Agrawal, 
{\em Nonlinear fiber optics}, 3rd ed.. Academic Press, San Diego, 2001.

\bibitem{agarwal2} G.P. Agrawal, {\em Applications of nonlinear fiber optics}, Academic Press, San Diego, 2001.

\bibitem{BCIR1}
O. Bang, P.L. Christiansen, F. If, K.O. 
Rasmussen, Y.B. Gaididei,  {\em Temperature effects in a
nonlinear model of monolayer Scheibe aggregates},
Phys. Rev. E, {\sl 49}, 4627--4636, (1994). 

\bibitem{BCIR2}
O. Bang, P.L. Christiansen, F. If, K.O. 
Rasmussen, Y.B. Gaididei,  {\em White Noise in the Two-dimensional
Nonlinear Schr\"odinger Equation}, Appl. Anal., {\sl 57}, 3--15, (1995).


\bibitem{cazenave} T. Cazenave,
{\em Semilinear Schr\"odinger equations},
Courant Lecture Notes in Mathematics, American Mathematical Society,
Courant Institute of Mathematical Sciences, 2003.

\bibitem{dBDL2}
A. de Bouard, A. Debussche, 
{\em A stochastic nonlinear Schr\"odinger equation with
multiplicative noise}, Comm. in Math. Phys., {\bf 205}, p.161--181 (1999).

\bibitem{dBDH1}
A. de Bouard, A. Debussche, 
{\em The stochastic nonlinear Schr\"odinger equation in $H^1$},
Stochastic Anal. Appl. {\bf 21}, p. 97--126 (2003).

\bibitem{dBDblowupadditif}
A. de Bouard, A. Debussche,
{\em   On the effect of a noise on the solutions of supercritical
     Schr\"odinger equation}, Prob. Theory and Rel. Fields, {\bf 123}, p. 76--96 
(2002).

\bibitem{dBDblowupmultip}
A. de Bouard, A. Debussche,
{\em Blow-up for the supercritical stochastic nonlinear
Schr\"odinger equation with multiplicative noise}, Ann. Probab., {\bf 33}, no. 3, 
p. 1078--1110 (2005).

\bibitem{dBD10} A. de Bouard, A. Debussche,
{\em The nonlinear Schrodinger equation with white noise dispersion}, Journal of Functional Analysis, 259, pp. 1300-1321 (2010).

\bibitem{dB-R} A. de Bouard, R. Fukuizumi, 
{\em Representation formula for stochastic Sch\"odinger evolution equations and applications},
Preprint.


\bibitem{DDM}
A. Debussche, L. DiMenza,
{\em Numerical simulation of focusing stochastic nonlinear
     Schr\"odinger equations}, Physica D, {\bf 162} (3-4), p. 131--154 (2002).

\bibitem{deb-gautier} A. Debussche, E. Gautier, 
{\em Small noise asymptotic of the timing jitter in soliton transmission}, Annals of Applied Probability 18, 178-208,  (2008).

\bibitem{durrett} R. Durrett, {\em Probablities: Theory and Examples}, Third Edition, Thomson, 2005.

\bibitem{FKLT}
G.E. Falkovich, I. Kolokolov, V. Lebedev. S.K. Turitsyn,
{\em Statistics of soliton-bearing systems with additive noise},
Phys. Rev. E, {\sl 63}, (2001).



\bibitem{Ga} J. Garnier,
{\em Stabilization of dispersion managed solitons in random optical fibers by strong
dispersion management}, 
Opt. Commun. {\bf 206}, p. 411--438 (2002).

\bibitem{gautier1} E. Gautier
{\em Large deviations and support results for the nonlinear Schrodinger equation with additive noise}, ESAIM: Probability and Statistics 9, 74-97,  (2005).

\bibitem{gautier2} E. Gautier
{\em Uniform large deviations for the nonlinear Schrodinger equation with multiplicative noise}, Stochastic Processes and Their Applications 115, 1904-1927 (2005).


\bibitem{Gi-Ve}
J. Ginibre, G. Velo,
{\em The global Cauchy problem for the nonlinear
Schr\"odinger equation revisited},
Ann. Inst. Henri Poincar\'e, Analyse Non Lin\'eaire
{\bf 2}, p. 309--327 (1985).

\bibitem{kato}
T. Kato, 
{\em On Nonlinear Schr\"odinger Equation}. Ann. Inst. H. Poincar\'e, Phys. Th\'eor.
{\sl 46}, p. 113--129 (1987).

\bibitem{konotop-vasquez}
{Konotop, V.; V\'azquez L. }
{\em Nonlinear random waves}; World Scientific Publishing Co., Inc.:
River Edge, N.J., 1994.



\bibitem{marty} 
R. Marty  {\em On a splitting scheme for the nonlinear Schr\"odinger equation in a random medium}, 
Commun. Math. Sci. {\bf 4}, no. 4, p. 679--705 (2006).

\bibitem{ozawa-tsutsumi} T. Ozawa, Y. Tsutsumi, 
{\em  Space-time estimates for null gauge forms and nonlinear Schr\"oinger equations},  Differential Integral Equations,  11,  no. 2, 201--222, (1998).

%
%

\bibitem{sulem-sulem} {C. Sulem, P.L. Sulem,}
{\em The Nonlinear Schr\"odinger Equation, Self-Focusing and
Wave Collapse};
Appl. Math. Sciences, Springer Verlag: New York, 1999.


\bibitem{tsutsumi}
Y. Tsutsumi, 
{\em $L^2$-solutions for nonlinear Schr\"odinger
equations and nonlinear groups}, 
Funk. Ekva. {\bf 30}, p. 115--125 (1987).

\bibitem{UK}
{Ueda, T.;  Kath W.L.}
{Dynamics of optical pulses in randomly birefrengent fibers}.
Physica D {\bf 1992}, {\sl 55}, 166--181.

\bibitem{ZGJT01} V. Zharnitsky, E. Grenier, C. Jones,  S. Turitsyn,
{\em Stabilizing effects of dispersion management},
Phys. D {\bf 152/153}, p. 794--817 (2001). 
\end{thebibliography}
\end{document}